\documentclass[11pt]{article}

\usepackage{geometry}
\usepackage[T1]{fontenc}
\usepackage[latin1]{inputenc}
\usepackage{setspace}
\usepackage{indentfirst}

\makeatletter

\usepackage{verbatim}

\makeatother

\usepackage[intlimits]{amsmath}
\usepackage{amssymb}
\usepackage{stmaryrd}
\usepackage{wasysym}
\usepackage{pifont}
\usepackage{graphicx}
\usepackage[all]{xy}

\usepackage{amsthm}

\theoremstyle{plain}
\newtheorem{prop}{Proposition}[section]
\newtheorem{lem}[prop]{Lemma}
\newtheorem{fact}[prop]{Fact}

\newtheorem{cor}[prop]{Corollary}
\newtheorem{them}[prop]{Theorem}

\newtheorem{question}[prop]{Question}

\newtheorem*{prop*}{Proposition}
\newtheorem*{lem*}{Lemma}
\newtheorem*{sublem*}{Sublemma}
\newtheorem*{cor*}{Corollary}
\newtheorem*{thm*}{Theorem}
\newtheorem*{hypo*}{Hypothesis}
\newtheorem*{question*}{Question}
\newtheorem*{conjecture*}{Conjecture}
\newtheorem*{scholum*}{Scholum}

\theoremstyle{definition}
\newtheorem{defn}[prop]{Definition}

\newtheorem{rmk}[prop]{Remark}

\newtheorem{ex}[prop]{Example}
\newtheorem{exs}[prop]{Examples}
\newtheorem*{defn*}{Definition}
\newtheorem*{con*}{Construction}
\newtheorem*{note*}{Note}
\newtheorem*{rmk*}{Remark}
\newtheorem*{rmks*}{Remarks}
\newtheorem*{ex*}{Example}
\newtheorem*{exs*}{Examples}

\theoremstyle{remark}

\newtheorem*{warning*}{Warning}
\newtheorem*{shortnote*}{Note}
\newtheorem*{claim*}{Claim}
\newtheorem*{axiom*}{Axiom}
\newtheorem*{example*}{Example}
\newtheorem*{examples*}{Examples}
\newtheorem*{remark*}{Remark}
\newtheorem*{remarks*}{Remarks}

\DeclareMathOperator{\ad}{ad}

\DeclareMathOperator{\idn}{id}

\DeclareMathOperator{\TTT}{T}
\DeclareMathOperator{\mult}{mult}

\renewcommand{\epsilon}{\varepsilon}
\newcommand{\rond}{\smalcirc}

\newcommand{\cty}{C^{\infty}}
\newcommand{\inv}{^{-1}}

\newcommand{\gradue}{^{\bullet}}

\newcommand{\sections}[1]{\Gamma(#1)}

\newcommand{\derlie}[1]{\mathcal{L}_{#1}}
\newcommand{\ip}[2]{\langle #1,#2 \rangle}
\newcommand{\lie}[2]{[#1,#2]}
\newcommand{\std}[2]{\llparenthesis #1,#2 \rrparenthesis}
\newcommand{\courant}[2]{\llbracket#1,#2\rrbracket}
\newcommand{\anchor}{\rho}
\newcommand{\derivation}{\delta}
\newcommand{\ttf}{\phi}
\newcommand{\ipJ}[2]{\langle #1,#2 \rangle_J}

\newcommand{\anchorJ}{\rho_J}

\newcommand{\lieN}[2]{[#1,#2]_N}

\newcommand{\DD}{\mathcal{D}}

\newcommand{\bivector}{\pi}
\newcommand{\bivectorN}{\pi_N}
\newcommand{\nijenhuis}{N}
\newcommand{\vf}{\mathfrak{X}}
\newcommand{\df}{\Omega}
\newcommand{\diese}{^{\sharp}}
\newcommand{\bemol}{_{\flat}}
\newcommand{\transpose}{^{\TTT}}
\newcommand{\interieur}[1]{i_{#1}}
\newcommand{\bmap}{B}
\newcommand{\cpn}{C^{\nijenhuis}_{\bivector\diese}}
\newcommand{\thalf}{\tfrac{1}{2}}

\newenvironment{rome}
{
\begin{enumerate}}
{\end{enumerate}
}

\usepackage{amscd}      

           {\nolinebreak $\Box$ \end{trivlist}}

\newcommand{\noprint}[1]{}

\renewcommand{\tilde}{\widetilde}

\newcommand{\toto}{\rightrightarrows}
\newcommand{\cinf}{\mathop{C^\infty}}

\newcommand{\XX}{{\mathfrak X}}

\newcommand{\cala}{{\cal A}}

\def\gpd{\rightrightarrows}

\newcommand{\ldiag}[1]%
       {\makebox[0cm]{${\scriptstyle#1}\downarrow\phantom{\scriptstyle#1}$}}
\newcommand{\ldiagup}[1]%
       {\makebox[0cm]{${\scriptstyle#1}\uparrow\phantom{\scriptstyle#1}$}}
\newcommand{\rdiag}[1]%
       {\makebox[0cm]{$\phantom{\scriptstyle#1}\downarrow{\scriptstyle#1}$}}
\newcommand{\sediagr}[1]%
       {\makebox[0cm]{$\phantom{\scriptstyle#1}\searrow{\scriptstyle#1}$}}
\newcommand{\nediagr}[1]%
       {\makebox[0cm]{$\phantom{\scriptstyle#1}\nearrow{\scriptstyle#1}$}}
\newcommand{\rdiagup}[1]%
       {\makebox[0cm]{$\phantom{\scriptstyle#1}\uparrow{\scriptstyle#1}$}}
\newcommand{\swdiag}[1]%
       {\makebox[0cm]{$\phantom{\scriptstyle#1}\swarrow{\scriptstyle#1}$}}
\newcommand{\sediag}[1]%
       {\makebox[0cm]{${\scriptstyle#1}\searrow\phantom{\scriptstyle#1}$}}
\newcommand{\nediag}[1]%
       {\makebox[0cm]{${\scriptstyle#1}\nearrow\phantom{\scriptstyle#1}$}}

\newcommand{\doublearrowstack}[2]%
                      {{{{\scriptstyle#1}\atop{\textstyle\longrightarrow}}\atop{{\textstyle\longrightarrow}\atop{\scriptstyle#2}}}}
\newcommand{\rightleftarrowstack}[2]%
                      {{{{\scriptstyle#1}\atop{\textstyle\longrightarrow}}\atop{{\textstyle\longleftarrow}\atop{\scriptstyle#2}}}}
\newcommand{\leftrightarrowstack}[2]%
                      {{{{\scriptstyle#1}\atop{\textstyle\longleftarrow}}\atop{{\textstyle\longrightarrow}\atop{\scriptstyle#2}}}}

\newcommand{\overtoparrow}%
{\makebox[0cm]{\beginpicture \setcoordinatesystem units
<.8cm,.4cm> point at 0 0 \setplotarea x from -3 to 3, y from 0 to
1 \setquadratic \plot -3 0 0 1 3 0 / \put{\vector(3,-1){0}}[Bl] at
3 0
\endpicture}}

\newcommand{\underbottomarrow}%
{\makebox[0cm]{\beginpicture \setcoordinatesystem units
<.8cm,.4cm> point at 0 0 \setplotarea x from -3 to 3, y from 0 to
1 \setquadratic \plot -3 1 0 0 3 1 / \put{\vector(3,1){0}}[Bl] at
3 1
\endpicture}}

\newcommand{\ses}[5]%
{0\longrightarrow#1\stackrel{#2}{ \longrightarrow}#3\stackrel{#4}{
\longrightarrow}#5\longrightarrow0}

\newcommand{\dt}[6]%
{#1\stackrel{#2}{longrightarrow}#3
\stackrel{#4}{\longrightarrow}#5 \stackrel{#6}{\longrightarrow}
#1[1]}

\newcommand{\cat}[1]%
{(\mbox{\rm #1})}

\newcommand{\gm}{\Gamma }

\newcommand{\be}{\begin{eqnarray*}}
\newcommand{\ee}{\end{eqnarray*}}
\newcommand{\smalcirc}{\mbox{\tiny{$\circ $}}}

\newcommand{\tpi}{{\tilde{\pi}}}
\newcommand{\tN}{{\tilde{N}}}
\newcommand{\find}{(TM\oplus T^*M,\ip{\cdot}{\cdot}_J,\std{\cdot}{\cdot}_J,\anchor_J)}
\newcommand{\tomega}{\tilde{\omega}}

\let\Vec=\overrightarrow
\let\ceV=\overleftarrow

\begin{document}

\title{Poisson Quasi-Nijenhuis Manifolds}
\author{\textsc{Mathieu Sti\'enon} \thanks{Francqui fellow of the Belgian American Educational Foundation} \\ Departement Mathematik \\ E.T.H. Zürich \\ Rämistrasse 101 \\ 8092 Zürich, Switzerland \\ \texttt{mathieu.stienon@math.ethz.ch} \and \textsc{Ping Xu} \thanks{Research supported by NSF grant DMS03-06665 and NSA grant 03G-142.} \\ Department of Mathematics \\ Pennsylvania State University \\ University Park, PA 16802\\USA \\ \texttt{ping@math.psu.edu} } 
\date{}
\maketitle

\begin{abstract}
We introduce the notion of Poisson quasi-Nijenhuis manifolds generalizing
the Poisson-Nijenhuis manifolds of Magri-Morosi.
We also investigate the integration problem
of Poisson quasi-Nijenhuis manifolds. In particular,
we prove that, under some topological assumption,
Poisson (quasi)-Nijenhuis manifolds are in one-one
correspondence with symplectic (quasi)-Nijenhuis groupoids.
As an application,  we study   generalized complex structures in
terms of Poisson quasi-Nijenhuis manifolds.
We prove that a generalized complex
manifold corresponds to a special class of Poisson quasi-Nijenhuis
structures.  As a consequence, we show that a generalized complex structure
integrates to a symplectic quasi-Nijenhuis groupoid recovering  a
theorem of Crainic.
\end{abstract}

\tableofcontents

\section{Introduction}

Poisson Nijenhuis structures were introduced 
by Magri and Morosi \cite{MR773513,MR791643} in their study of
bi-Hamiltonian systems, and intensively studied
by many authors \cite{MR1077465,MR1390832}.
Recall that a Poisson Nijenhuis manifold  consists of 
a triple $(M,\pi,N)$, where $M$ is a manifold endowed with 
a Poisson bivector field $\pi$, and a $(1,1)$-tensor 
$N$ whose Nijenhuis torsion vanishes, i.e.
 $$[NX,NY]-N([NX,Y]+[X,NY]-N[X,Y])=0, \ \ \ \forall X, Y\in \mathfrak{X}(M),$$  
together with some compatibility condition between $\pi$ and $N$. Poisson 
Nijenhuis structures are very important in the study of integrable systems 
since they produce bi-Hamiltonian systems \cite{MR773513,MR1077465}. 

As observed by Kosmann-Schwarzbach \cite{MR1421686}, 
given a Poisson Nijenhuis manifold $(M,\pi,N)$, $((T^*M)_\pi,(TM)_N)$ 
constitutes a Lie bialgebroid, where 
$(T^*M)_\pi$ is equipped with the standard cotangent Lie algebroid 
structure induced by the Poisson tensor $\pi$ while $(TM)_N$ 
is the deformed Lie algebroid on $TM$ induced by the Nijenhuis endomorphism $N$.
 Indeed it is proved in \cite{MR1421686} that the Lie bialgebroid 
condition on $((T^*M)_\pi,(TM)_N)$ is equivalent
 to the triple $(M,\pi,N)$ being Poisson Nijenhuis.

The main goal of the present paper is
to  introduce the notion of Poisson quasi-Nijenhuis 
structures. By definition, a Poisson quasi-Nijenhuis manifold 
is a quadruple $(M,\pi,N,\ttf)$, where $M$ is manifold endowed 
with a Poisson bivector field $\pi$, a $(1,1)$-tensor
$N$ and a closed 3-form $\ttf$ such that
$\pi$ and $N$ are compatible (in the usual Poisson-Nijenhuis sense) and 
$$[NX,NY]-N([NX,Y]+[X,NY]-N[X,Y])=\pi^\sharp(i_{X\wedge Y}\ttf), 
\quad\forall X,Y\in\mathfrak{X}(M) .$$

Recall that Lie bialgebroids 
are pairs of transverse Dirac structures in
 a Courant algebroid \cite{MR1472888}. 
When one of the two maximal isotropic direct summands fails
 to be Courant involutive, 
this becomes a quasi-Lie bialgebroid \cite{MR1936572,math.DG/9910078}. Alternatively, 
a quasi-Lie bialgebroid is equivalent to the following data: a Lie 
algebroid $A$ together with a degree 1 derivation $\delta$ 
of the associated Gerstenhaber algebra 
$\big(\Gamma(\wedge^\bullet A),\wedge,[\cdot,\cdot]\big)$ such that 
$\delta^2=[\ttf,\cdot]$ and $\delta\ttf=0$ for some $\ttf\in\Gamma(\wedge^3 A)$ \cite{math.DG/0507396}.
We prove

\theoremstyle{plain}
\newtheorem*{thmA}{Theorem A}
\begin{thmA}
Given $(M,\pi,N,\ttf)$, the following are equivalent
\begin{itemize}
\item $(M,\pi,N,\ttf)$ is a Poisson quasi-Nijenhuis  manifold;
\item $\big((T^*M)_{\pi},(TM)_N,\ttf\big)$ is a quasi-Lie bialgebroid.
\end{itemize}
\end{thmA}

It is well known that the global object
corresponding to a Poisson manifold is a symplectic
groupoid \cite{MR996653,MR866024}. It is natural
to ask what is the global object integrating
a Poisson Nijenhuis manifold. We prove

\theoremstyle{plain}
\newtheorem*{thmC}{Theorem B}
\begin{thmC}
The base manifold of a  symplectic Nijenhuis groupoid
is a Poisson Nijenhuis manifold.
Moreover, 
there is a
one-one correspondence between $t$-connected and $t$-simply
connected symplectic Nijenhuis groupoids
$(\gm\toto M,\tomega,\tN)$ and  integrable
Poisson Nijenhuis manifolds $(M, \pi, N)$.
\end{thmC}

By a symplectic Nijenhuis groupoid, we mean a symplectic groupoid
$(\gm\toto M,\tomega)$ equipped with a multiplicative
$(1,1)$-tensor $\tN:T\gm\to T\gm$ such that
$(\gm,\tomega,\tN)$ is a symplectic Nijenhuis structure.
The main idea of the proof of Theorem B can be outlined as
follows. One proves that Poisson Nijenhuis
structures on a manifold $M$ are in one-one
correspondence with Lie bialgebroids $((T^*M)_\pi,\delta)$
satisfying the condition  that $[\delta,d]=0$, where
$d$ is the de Rham differential on $M$.
The latter are the infinitesimal of symplectic Nijenhuis groupoids,
as can be shown using the universal lifting theorem \cite{math.DG/0507396}.

The same method can be used to prove an analogous result
for Poisson quasi-Nijenhuis manifolds.

\theoremstyle{plain}
\newtheorem*{thmD}{Theorem C}
\begin{thmD}
The base manifold of a symplectic quasi-Nijenhuis groupoid
is a Poisson quasi-Nijenhuis manifold.
Moreover
there is a
one-one correspondence between  $t$-connected and $t$-simply
connected symplectic quasi-Nijenhuis groupoids
$\big(\gm\toto M,\tomega,\tN,t^*\ttf-s^*\ttf\big)$
and integrable Poisson quasi-Nijenhuis manifolds $(M,\pi,N,\ttf)$.
\end{thmD}

A symplectic quasi-Nijenhuis groupoid is a symplectic groupoid
$(\gm\toto M,\tomega)$ equipped with a  multiplicative $(1,1)$-tensor
$\tN:T\gm\to T\gm$ and a closed $3$-form $\ttf\in\Omega^3(M)$
such that $\big(\gm,\tomega,\tN,t^*\ttf-s^*\ttf\big)$ is a symplectic
quasi-Nijenhuis structure.

As an application, we study generalized complex structures
in terms of Poisson quasi-Nijenhuis structures.
The notion of generalized complex structures was
 introduced by Hitchin \cite{MR2013140} and studied by Gualtieri \cite{math.DG/0401221} motivated by the study of mirror symmetry.
 It comprises both symplectic and complex structures as extreme cases.
We show that  on a generalized complex
manifold $(M,J)$, where
$$J=\begin{pmatrix} N & \pi^{\sharp} \\ \sigma_{\flat} & -N^* \end{pmatrix}$$ with $N^2+\pi^{\sharp} \sigma_{\flat }=-\idn$,
 the building units $\pi$, $N$ and $\sigma$ of $J$ do exactly 
determine a Poisson quasi-Nijenhuis structure. Indeed, the endomorphism $N$ can be used to define a derivation $d_N$ of the Gerstenhaber algebra
associated to the Lie algebroid  $(T^*M)_\pi$. We prove

\theoremstyle{plain}
\newtheorem*{thmB}{Theorem D}
\begin{thmB}
The following are equivalent
\begin{itemize}
\item $J$ is a generalized complex structure;
\item $(M,\pi,N,d\sigma)$ is a Poisson quasi-Nijenhuis structure such that
$$(TM)_N\oplus (T^*M)_{\pi}\xrightarrow{J}TM\oplus T^*M$$
is a Courant algebroid isomorphism.
\end{itemize}
\end{thmB}

A similar result (in a different form) was already proved
 by Crainic using a direct argument \cite{math.DG/0412097}.

Since a generalized complex structure corresponds to a 
quasi-Nijenhuis manifold according to Theorem D,
as a consequence, we prove 

\theoremstyle{plain}
\newtheorem*{thmE}{Theorem E}
\begin{thmE}
Let $J$ be a generalized complex structure as given by
Eq. \eqref{eq:J}, and $(\gm\toto M,\tomega)$ a $t$-connected and
$t$-simply connected symplectic groupoid integrating
$(T^*M)_\pi$. Then there is a multiplicative
$(1,1)$-tensor $\tN$ on $\gm$ such that
$\big(\gm\toto M,\tomega,\tN,t^*d\sigma-s^*d\sigma\big)$ is a
symplectic quasi-Nijenhuis groupoid.
\end{thmE}

This result, in a disguised form, was also proved by Crainic \cite{math.DG/0412097}
using a different method.

\paragraph{Notations}
We denote the bracket on the sections of a Courant algebroid by $\courant{\cdot}{\cdot}$, except for the standard Courant bracket on $TM\oplus T^*M$, which is denoted by $\std{\cdot}{\cdot}$. The Lie bracket of vector fields and its extension to polyvector fields (i.e. the Schouten bracket) are denoted by $\lie{\cdot}{\cdot}$. Any bundle map $B:T^*M\to TM$ induces a bracket on the space of 1-forms (see Eq. \eqref{eq:pi}). It is denoted by $\lie{\cdot}{\cdot}_{B}$ as well as its extension to the space of differential forms of all degrees. Finally, if $\courant{\cdot}{\cdot}$ is a bracket on the space of sections of a vector bundle $E$ of which $J$ is a bundle endomorphism, then its deformation by $J$ is denoted by $\courant{\cdot}{\cdot}_J$ (see Eq. \eqref{eq:N2}).

\paragraph{Acknowledgments}  
We  would like to thank several institutions
for their hospitality while work on this project was being done:
Erwin Schroedinger  International Institute for Mathematical Physics
(Sti\'enon and Xu), and Universit\'e Pierre et Marie Curie (Xu).
Sti\'enon is grateful to the \'Emile Francqui and Belgian American 
Educational Foundations for supporting his stay at the Pennsylvania 
State University where this work was completed in 2005.
 We would also like to thank Marius Crainic and
 Jim Stasheff for many useful discussions, and the referee for helpful
suggestions to improve the presentation of the paper.

\section{Preliminaries}

\begin{defn}[\cite{MR1472888}]
A Courant algebroid is a triple consisting of 
\begin{itemize}
\item a vector bundle $E\to M$ equipped with a non degenerate symmetric bilinear form $\ip{\cdot}{\cdot}$,
\item a skew-symmetric bracket $\courant{\cdot}{\cdot}$ on $\sections{E}$, and
\item a smooth bundle map $E\xrightarrow{\anchor}M$ called the anchor, which
induces a natural differential operator $\DD :\cty(M)\to\sections{E}$
defined by
$$\ip{\DD f}{A}=\tfrac{1}{2}\anchor(A)f$$ 
for all $f \in \cty(M)$ and $A\in\sections{E}$.
\end{itemize}
These structures must be compatible in the following sense: $\forall A,B,C\in\sections{E}$ and $\forall f,g\in \cty(M)$, 
\begin{itemize}
\item $\anchor(\courant{A}{B})=\lie{\anchor(A)}{\anchor(B)}$,
\item $\courant{\courant{A}{B}}{C}+\courant{\courant{B}{C}}{A}+\courant{\courant{C}{A}}{B}=
\tfrac{1}{3}\DD\big(\ip{\courant{A}{B}}{C}+\ip{\courant{B}{C}}{A}+\ip{\courant{C}{A}}{B}\big)$, 
\item $\courant{A}{fB}=f\courant{A}{B}+\big(\anchor(A)f\big)B-\ip{A}{B}\DD f$,
\item $\anchor\rond\DD=0$, i.e. $\ip{\DD f}{\DD g}=0$, 
\item $\anchor(A)\ip{B}{C}=\ip{\courant{A}{B}+\DD \ip{A}{B}}{C}+\ip{B}{\courant{A}{C}+\DD \ip{A}{C}}$.
\end{itemize}
Note that a Courant algebroid is not a Lie algebroid as the Jacobi identity is not satisfied.
\end{defn}

\begin{ex}[\cite{MR998124}]
\label{ex:1.2}
The generalized tangent bundle $TM\oplus T^*M$ of a manifold $M$ is
a Courant algebroid,
where the anchor is the
projection onto the first component and the pairing
and bracket are given, respectively, by
\begin{gather}
\ip{X+\xi}{Y+\eta}=\tfrac{1}{2}\big(\xi(Y)+\eta(X)\big), \label{eq:pairing} \\
\std{X+\xi}{Y+\eta}=\lie{X}{Y}+\derlie{X}\eta-\derlie{Y}\xi+\tfrac{1}{2}\big(\xi(Y)-\eta(X)\big), \label{eq:bracket}
\end{gather}
$\forall X,Y \in \vf(M)$, $\forall \xi,\eta \in \df^1(M)$.
\end{ex}

\begin{defn}
A Dirac structure is a smooth subbundle $L$ of a Courant algebroid $E$, which is maximal isotropic with respect to $\ip{\cdot}{\cdot}$ and whose space of sections $\sections{L}$ is closed under $\courant{\cdot}{\cdot}$. It is thus naturally a Lie algebroid.
\end{defn}

It is well-known \cite{MR1675117} that a Lie algebroid $(A,\lie{\cdot}{\cdot}_A,\anchor_A)$
 gives rise to a Gerstenhaber algebra $(\sections{\wedge\gradue A},\wedge,\lie{\cdot}{\cdot}_A)$,
 and a degree 1 derivation $\derivation_{A}$ of the graded commutative algebra $(\sections{\wedge\gradue A^*},\wedge)$ such that $(\derivation_{A})^2=0$. Here $\derivation_{A}$ is given by
\begin{multline}
\label{eq:deltaA}
(\derivation_{A}\alpha)(X_0,X_1,\cdots,X_n)=\sum_{i=0}^n (-1)^i (\anchor_A X_i)
 \alpha(X_0,\cdots,\widehat{X_i},\cdots,X_n) \\
+ \sum_{i<j} (-1)^{i+j} \alpha(\lie{X_i}{X_j}_A,X_0,\cdots,\widehat{X_i},\cdots,\widehat{X_j},\cdots,X_n).
\end{multline}

A Lie bialgebroid \cite{MR1746902,MR1262213} is a pair of Lie algebroid structures on $A$ and its dual $A^*$ such that $\derivation_{A^*}$ is a derivation of the Gerstenhaber algebra $(\sections{\wedge\gradue A},\wedge,\lie{\cdot}{\cdot}_A)$ or, equivalently, such that $\derivation_{A}$ is a derivation of the Gerstenhaber algebra $(\sections{\wedge\gradue A^*},\wedge,\lie{\cdot}{\cdot}_{A^*})$.
Since the bracket $\lie{\cdot}{\cdot}_{A^*}$ can be recovered from the derivation $\derivation_{A^*}$, one is led to the following alternative definition.

\begin{defn}
A Lie bialgebroid is a pair $(A,\derivation)$ consisting
 of a Lie algebroid $(A,\lie{\cdot}{\cdot}_A,\anchor_A)$ and a degree 1 derivation $\derivation$ of the Gerstenhaber algebra $(\sections{\wedge\gradue A},\wedge,\lie{\cdot}{\cdot}_A)$ such that $\derivation^2=0$.
\end{defn}

More generally, we can speak about quasi-Lie bialgebroids \cite{MR1936572,math.DG/0507396}.

\begin{defn}[\cite{math.DG/0507396}]
A quasi Lie bialgebroid is a triple $(A,\derivation,\ttf)$ 
consisting of a Lie algebroid $A$, a degree 1 derivation $\derivation$ 
of the Gerstenhaber algebra $(\sections{\wedge\gradue A},\wedge,\lie{\cdot}{\cdot}_A)$ 
and an element $\ttf\in\sections{\wedge^3 A}$ such that 
$\derivation^2=\lie{\ttf}{\cdot}_A$ and $\derivation\ttf=0$.
\end{defn}

The link between Courant, Lie bi- and quasi Lie bialgebroids is given by the following 

\begin{them}[\cite{MR1472888,MR1936572,math.DG/9910078}]
\label{thm:1.6}
\begin{enumerate}
\item 
There is a 1-1 correspondence between Lie bialgebroids 
and pairs of transversal Dirac structures in a Courant algebroid.
\item There is a 1-1 correspondence between quasi Lie bialgebroids 
and  Dirac structures with transversal isotropic complements
in a Courant algebroid.
\end{enumerate}
\end{them}

\begin{proof}
The proof of (i) can be found in \cite{MR1472888}, and 
(ii) was proved in  \cite{MR1936572,math.DG/9910078}. Below
we give an explicit formula describing such
a correspondence, which will be needed later.

Let $(A,\derivation,\ttf)$ be a quasi Lie bialgebroid.
Let $\rho_{A*} :A^*\to TM$ be the bundle map given by
$$\rho_{A^*} (\xi ) (f)=\xi(\derivation f), \quad \forall \xi\in
A^*, \ \forall f\in C^{\infty}(M).$$
Introduce a bracket on $\sections{A^*}$ by
$$[\xi, \eta ]_{A^*}(X) =  (\rho_{A^*}\xi)(\eta X) -(\rho_{A^*}\eta)(\xi X)
-(\delta X)(\xi, \eta).$$
Note that $(A^*, \rho_{A^*}, [\cdot, \cdot ]_{A*})$ is
in general not a Lie algebroid.
Let $E= A^*\oplus A$ and  $\rho:E\to TM$ be the bundle map
$$\rho (\xi  +X )= \rho_{A^*} (\xi )+\rho_A (X) .$$
Define a non-degenerate symmetric pairing on $E$ by
$$\ip{\xi+X}{\eta +Y}=\tfrac{1}{2}\big(\xi(Y)+\eta(X)\big),$$
and a bracket $\courant{\cdot}{\cdot}$ on $\sections{E}$ by
\begin{gather}
\courant{X}{Y} = [X,Y]_A, \nonumber \\
\courant{\xi}{\eta} = [\xi, \eta]_{A^*} +\ttf (\xi, \eta, \cdot ), \nonumber \\
\courant{X}{\xi} = \big(i_X\derivation_{A^*}\xi+\thalf\derivation_{A^*}(\xi X)\big)
-\big(i_{\xi}\derivation_A X+\thalf\derivation_A(\xi X)\big), \label{eq:5quater}
\end{gather}
for all $X,Y\in\sections{A}$ and $\xi,\eta\in\sections{A^*}$.
Here $\derivation_{A*}:\sections{\wedge^{\bullet}A^*}\to\sections{\wedge^{\bullet+1}A^*}$
is the derivation given by Eq. \eqref{eq:deltaA}.
Then $(E,  \ip{\cdot}{\cdot},\courant{\cdot}{\cdot},\anchor)$ is a 
Courant algebroid.

Conversely, assume that  $(E,\ip{\cdot}{\cdot},\courant{\cdot}{\cdot},\anchor)$ 
is a Courant algebroid, and $A$ is a Dirac structure with an isotropic 
complement $B$. The duality pairing $$A\otimes B\to\mathbb{R}:X\otimes\xi\mapsto 2\ip{\xi}{X}$$ identifies $B$ with $A^*$. 
Let $\phi$ be the element in  $\sections{\wedge^3 A}$ defined by 
\begin{equation} \ttf (\xi, \eta, \zeta)= 2 \ip{\courant{\xi}{\eta}}{\zeta}, \quad \forall
\xi, \eta, \zeta \in \sections{B},
\label{eq:5ter} \end{equation}
$\rho_B=\rho|_B$ be the restriction of $\anchor$ to 
$B$ and $[\cdot,\cdot]_B$ be the bracket on $\sections{B}$ such that 
\begin{equation} 
\courant{\xi}{\eta}-[\xi,\eta]_B\in\sections{A},\quad\forall \xi,\eta\in\sections{B}.
\label{eq:5bis} \end{equation} 

Define a derivation $\derivation:\sections{\wedge^{\bullet}A} 
(\cong\sections{\wedge^{\bullet}B^*})\to\sections{\wedge^{\bullet+1} A}
(\cong\sections{\wedge^{\bullet+1}B^*}$ as in Eq. \eqref{eq:deltaA}. The
triple $(A,\derivation,\ttf)$ becomes a quasi Lie bialgebroid.
\end{proof}

\section{Poisson quasi-Nijenhuis manifolds}

Let $M$ be a smooth manifold, $\bivector$ a Poisson bivector field, 
and $\nijenhuis:TM\to TM$ a $(1,1)$-tensor.

\begin{defn}[\cite{MR1421686}]
The bivector field $\bivector$ and the tensor $N$ are said to be compatible \cite{MR1077465} if 
\begin{equation} 
N\smalcirc\pi\diese=\pi\diese\smalcirc N\transpose\qquad\text{and}\qquad\cpn=0, \label{eq:newnir}
\end{equation}
where 
$$\cpn(\alpha,\beta):=\lie{\alpha}{\beta}_{N\bivector\diese}-\big(\lie{N\transpose \alpha}{\beta}_{\bivector\diese}
+\lie{\alpha}{N\transpose\beta}_{\bivector\diese}-N\transpose\lie{\alpha}{\beta}_{\bivector\diese}\big)$$
and 
\begin{equation} \label{eq:pi}
\lie{\alpha}{\beta}_{\bmap}:=\derlie{\bmap\alpha}(\beta)-\derlie{\bmap\beta}(\alpha)
-d\big(\beta(\bmap\alpha)\big)
\end{equation}
for all $\alpha,\beta\in\df^1(M)$ and any bundle map $\bmap:T^*M\to TM$.
\end{defn}

The $(1,1)$-tensor $N$ is said to have zero Nijenhuis torsion if 
$$\lie{\nijenhuis X}{\nijenhuis Y}-\nijenhuis\big(\lie{\nijenhuis X}{Y}
+\lie{X}{\nijenhuis Y}-\nijenhuis\lie{X}{Y}\big)=0,\quad\forall X,Y\in\vf({M}).$$

In \cite{MR900387}, Magri and Morosi defined Poisson Nijenhuis manifolds 
as triples $(M,\bivector,\nijenhuis)$ such that $\pi$ and $N$ are compatible 
and the Nijenhuis torsion of $N$ vanishes.

This definition is motivated by the following 

\begin{fact}[\cite{MR1077465,MR1390832}] \label{facts}
Assume that $\pi\in\XX^2(M)$ is a Poisson tensor 
and $N:TM\to TM$ a $(1,1)$-tensor on $M$.
The tensor $\pi_N$ defined by 
$$\pi_N(\alpha,\beta):=\beta(N\pi\diese\alpha),\quad\forall\alpha,\beta\in\Omega^1(M)$$
is skew-symmetric if, and only if, $N\smalcirc \pi\diese=\pi \smalcirc \diese N^T$.
In this case, we have
\begin{enumerate}
\item $\lie{\pi}{\pi_N}=0$ if $C^N_{\pi\diese}=0$, and the converse
 is true if $\pi$ is non-degenerate;
\item $\lie{\pi_N}{\pi_N}=0$ if, and only if, the Nijenhuis torsion of $N$ vanishes.
\end{enumerate}
\end{fact}

Hence, any Poisson Nijenhuis manifold $(M,\bivector,N)$ is endowed with a 
bi-Hamiltonian structure $(\pi, \pi_N )$, i.e. 
$$\lie{\bivector}{\bivector}=0,
\qquad\lie{\bivector}{\bivectorN}=0,
\qquad\lie{\bivectorN}{\bivectorN}=0.$$

Similarly, one can define Poisson quasi-Nijenhuis manifolds.

Let $i_N$ be the degree 0 derivation of $(\Omega^{\bullet}(M),\wedge)$ defined by
$$(i_N\alpha)(X_1,\cdots,X_p)=\sum_{i=1}^p\alpha(X_1,\cdots,NX_i,\cdots,X_p), \quad\forall\alpha\in\Omega^p(M).$$

\begin{defn} \label{def:42}
A Poisson quasi-Nijenhuis manifold is a quadruple $(M,\bivector,\nijenhuis,\ttf)$, 
where $\bivector\in\vf^2(M)$ is a Poisson bivector field, $N:TM\to TM$ is a $(1,1)$-tensor 
compatible with $\bivector$, and $\ttf$ is a closed 3-form on $M$ such that
$$\lie{\nijenhuis X}{\nijenhuis Y}-\nijenhuis\big(\lie{\nijenhuis X}{Y}+\lie{X}{\nijenhuis Y}
-\nijenhuis\lie{X}{Y}\big)=\bivector\diese(\interieur{X\wedge Y}\ttf), \quad\forall X,Y\in\vf(M)$$ 
and $i_N\ttf$ is closed.
\end{defn}

It is well known that, on a Poisson manifold $(M,\bivector)$, the bracket on $\df^1(M)$ associated to the bundle map $\pi\diese$ through Eq. \eqref{eq:pi} makes $T^*M$ into a Lie algebroid with anchor $\bivector\diese:T^*M\to TM$. 
The usual cotangent bundle will be denoted by $(T^*M)_\pi$ when equipped with this Lie algebroid structure.
More precisely, we have the following 

\begin{fact}[\cite{MR996653}] \label{fact}
Let $\pi$ be a bivector field on $M$. Then $\lie{\pi}{\pi}=0$ if, and only if, $(T^*M)_{\pi}$ is a Lie algebroid.
\end{fact}

On the other hand, defining a bracket 
$\lieN{\cdot}{\cdot}$ on $\XX(M)$ by $$[X,Y]_N =[NX,Y]+[X,NY]-N[X,Y],\quad\forall X,Y\in\XX(M)$$
as in \cite{MR1421686}, and considering $N:TM\to TM$ as an anchor map, we obtain a degree 1 derivation 
$d_N$ of $(\Omega^{\bullet}(M),\wedge)$ inspired by Eq. \eqref{eq:deltaA}:
\begin{multline}
(d_N\alpha)(X_0,X_1,\cdots,X_n)=\sum_{i=0}^n (-1)^i (N X_i)
 \alpha(X_0,\cdots,\widehat{X_i},\cdots,X_n) \\
+ \sum_{i<j} (-1)^{i+j} \alpha(\lie{X_i}{X_j}_N,X_0,\cdots,\widehat{X_i},\cdots,\widehat{X_j},\cdots,X_n).
\label{eq:8bis} \end{multline}

Moreover, as proved in \cite{MR1421686}, we have the following identity 
\begin{equation} d_N=[i_N, d]=i_N\smalcirc d-d \smalcirc i_N \label{eq:tgif} .\end{equation}

The following proposition extends a result of Kosmann-Schwarzbach \cite[Proposition 3.2]{MR1421686}.

\begin{prop}
\label{prop:qN}
The quadruple $(M,\bivector,\nijenhuis,\ttf)$ is a 
Poisson quasi-Nijenhuis manifold if, and only if,  
$\big((T^*M)_\pi, d_N,\ttf\big)$ is a quasi Lie bialgebroid 
and $\phi$ is a closed $3$-form.
\end{prop}

This is an immediate consequence of  Fact \ref{fact} and the following two lemmas.

\begin{lem}[{\cite[Proposition 3.2]{MR1421686}}] \label{lem:4}
Assume that $\pi\in\vf^2(M)$ is a Poisson tensor and $N:TM\to TM$ a $(1,1)$-tensor on $M$.
The differential $d_N$ is a derivation of the graded Lie algebra $(\Omega^{\bullet}(M),\lie{\cdot}{\cdot}_{\pi\diese})$ if, and only if, $\pi$ and $N$ are compatible.
\end{lem}

\begin{lem} \label{lem:5}
Let $(M,\pi)$ be a Poisson manifold and $N:TM\to TM$ a $(1,1)$-tensor compatible with $\pi\diese$. Then $d_N^2=\lie{\phi}{\cdot}_{\pi\diese}$ if, and only if, $$\lie{NX}{NY}-N(\lie{NX}{Y}+\lie{X}{NY}-N\lie{X}{Y})=\bivector\diese (i_{X\wedge Y}\phi), \quad\forall X,Y\in\vf(M)$$
and $\pi^{\#} \smalcirc (d\ttf)_\flat=0$, where $(d\ttf)_\flat: \wedge^3 TM\to T^*M$
is the bundle map defined by $(d\ttf)_\flat (u, v, w)=i_{u\wedge v\wedge w} d\ttf$, 
$\forall u, v, w\in TM$. 
\end{lem}

\begin{proof}
It follows from an easy computation that 
$$\big(d_N^2 f-\lie{\ttf}{f}_{\pi\diese}\big)(X,Y)=
(df)\big(\lie{NX}{NY}-N(\lie{NX}{Y}+\lie{X}{NY}-N\lie{X}{Y})-
\bivector\diese (i_{X\wedge Y}\phi)\big)$$
for all $f\in\cty(M)$. Moreover, since $d\rond d_N+d_N\rond d=0$, one has 
$$d_N^2 (df)-\lie{\ttf}{df}_{\pi\diese}=d(d_N^2 f)-\big(d\lie{\ttf}{f}_{\pi\diese}-\lie{d\ttf}{f}_{\pi\diese}\big)
= d(d_N^2 f-\lie{\ttf}{f}_{\pi\diese})+\lie{d\ttf}{f}_{\pi\diese} .$$
Hence, $d_N^2-\lie{\ttf}{\cdot}_{\pi\diese}$ vanishes on 0- and exact 1-forms if, and only if, 
$$\lie{NX}{NY}-N(\lie{NX}{Y}+\lie{X}{NY}-N\lie{X}{Y})=\bivector\diese (i_{X\wedge Y}\phi),\quad\forall X,Y\in\vf(M)$$ and 
$\lie{d\ttf}{f}_{\pi\diese}=0, \ \forall f\in C^\infty (M)$.
The latter is easily seen to be equivalent to
$\pi^{\#}\smalcirc (d\ttf)_\flat=0$. 
And in this case, since both $d_N^2$ and $\lie{\ttf}{\cdot}_{\pi\diese}$ are derivations with respect to $\wedge$, we get $d_N^2=\lie{\ttf}{\cdot}_{\pi\diese}$.
\end{proof}

As an immediate consequence, we obtain the following
result of Kosmann-Schwarzbach \cite{MR1421686}.
\begin{cor}
\label{cor:yvette}
The triple $(M,\bivector,\nijenhuis)$ is a Poisson Nijenhuis manifold 
if, and only if,  $((T^*M)_\pi, d_N)$ is a  Lie bialgebroid.
\end{cor}

We now turn our attention to the particular case where the Poisson bivector field $\pi$ is 
non-degenerate. Together with Lemma \ref{lem:4}, the  following two lemmas give another proof of the equivalence between the relation $\lie{\pi}{\pi_N}=0$ and the compatibility condition \eqref{eq:newnir} when $\pi$ is non-degenerate (see Fact \ref{facts}).

\begin{lem} 
\label{lem:PN}
Assume that $\pi\in\XX^2(M)$ is a Poisson tensor and $N: TM\to TM$ 
a $(1,1)$-tensor on $M$. Then $\pi_N$ is a bivector field such that $[\pi,\pi_N]=0$ if, and only if, all the squares in the following diagram commute.
\begin{equation} \xymatrix{ 0 \ar[r] & \cty(M) \ar[d]^{\idn} \ar[r]^{d_N} & \Omega^1(M) \ar[r]^{d_N} \ar[d]^{\bivector\diese} & 
\Omega^2(M) \ar[r]^{d_N} \ar[d]^{\bivector\diese} & \Omega^3(M) \ar[r]^{d_N} \ar[d]^{\bivector\diese} & \cdots \\ 
0 \ar[r] & \cty(M) \ar[r]_{\lie{\pi_N}{\cdot}} & \vf^1(M) \ar[r]_{\lie{\pi_N}{\cdot}} & 
\vf^2(M) \ar[r]_{\lie{\pi_N}{\cdot}} & \vf^3(M) \ar[r]_{\lie{\pi_N}{\cdot}} & \cdots } \label{eq:comdiag} \end{equation}
\end{lem}

\begin{proof}
We have $\pi\diese N^T=N \pi\diese$ (i.e. $\pi_N$ is a bivector field) if, and only if, 
$\forall f\in\cty(M)$,
\begin{align}
& \pi\diese N^T df= N\pi\diese df \nonumber \\
\Leftrightarrow \quad & \pi\diese i_N df = \pi_N\diese df \nonumber \\ 
\Leftrightarrow \quad & \pi\diese d_N f = \lie{\pi_N}{f} . \label{eq:a}
\end{align}
And $\lie{\pi_N}{\pi}=0$ is equivalent to
\begin{align}
& \lie{\pi_N}{\pi}\diese (df)=0 \nonumber \\
\Leftrightarrow \quad & \lie{\lie{\pi_N}{\pi}}{f}=0 \nonumber \\
\Leftrightarrow \quad & \lie{\lie{\pi_N}{f}}{\pi}+\lie{\pi_N}{\lie{\pi}{f}}=0 \nonumber \\
\Leftrightarrow \quad & \lie{\pi_N\diese df}{\pi}+\lie{\pi_N}{\pi\diese df}=0 \nonumber \\
\Leftrightarrow \quad & \lie{\pi}{\pi\diese N^T df}=\lie{\pi_N}{\pi\diese df} \nonumber \\
\Leftrightarrow \quad & \lie{\pi}{\pi\diese (i_N df)}=\lie{\pi_N}{\pi\diese df} \nonumber \\
\Leftrightarrow \quad & \pi\diese d(i_N df)=\lie{\pi_N}{\pi\diese df} \nonumber \\
\Leftrightarrow \quad & \pi\diese d_N (df)=\lie{\pi_N}{\pi\diese df} \label{eq:b} 
\end{align}
for all $f\in\cinf(M)$.
Since both $\pi\diese\rond d_N$ and $\lie{\pi_N}{\pi\diese(\cdot)}$ are derivations of $(\Omega^{\bullet}(M),\wedge)$, the equivalence follows from Eqs. \eqref{eq:a}-\eqref{eq:b}.
\end{proof}

\begin{lem} \label{lem:3}
Assume that $\pi\in\vf^2(M)$ is a non-degenerate Poisson tensor, and $N:TM\to TM$ is a $(1,1)$-tensor on $M$. If $\pi_N$ is a bivector field and Diagram \eqref{eq:comdiag} commutes, then $d_N$ is a derivation of $\lie{\cdot}{\cdot}_{\pi\diese}$.
\end{lem}

\begin{proof}
Since $\pi$ is Poisson, we have $$\pi\diese\lie{\alpha}{\beta}_{\pi\diese}=\lie{\pi\diese\alpha}{\pi\diese\beta}, \quad \forall \alpha,\beta\in\Omega^{\bullet}(M) .$$
Then, the Jacobi identity for the Schouten bracket gives 
$$\lie{\pi_N}{\pi\diese\lie{\alpha}{\beta}_{\pi\diese}}=\lie{\lie{\pi_N}{\pi\diese\alpha}}{\pi\diese\beta}
+\lie{\pi\diese\alpha}{\lie{\pi_N}{\pi\diese\beta}} ,$$ which can be rewritten as 
$$\pi\diese d_N\big(\lie{\alpha}{\beta}_{\pi\diese}\big)=\pi\diese\big(\lie{d_N\alpha}{\beta}_{\pi\diese}+\lie{\alpha}{d_N\beta}_{\pi\diese}\big)$$ since $\pi\diese\rond d_N=\lie{\pi_N}{\pi\diese(\cdot)}$. 
The conclusion follows from the invertibility of $\pi\diese$.
\end{proof}

The previous lemmas are used to prove the following

\begin{prop} \label{prop:3.10}
\begin{enumerate}
\item Let $(M,\pi,N,\ttf)$ be a Poisson quasi-Nijenhuis manifold. Then, 
\begin{equation} \lie{\pi}{\pi_N}=0, \label{eq:pqn2} \end{equation} and
\begin{equation} \lie{\pi_N}{\pi_N}=2 \pi\diese(\ttf) 
. \label{eq:pqn1} \end{equation}
\item Conversely, assume that $\pi\in\vf^2(M)$ is a non-degenerate Poisson bivector field, $N:TM\to TM$ is a $(1,1)$-tensor and $\ttf$ is a closed 3-form. If
 Eqs. \eqref{eq:pqn2}-\eqref{eq:pqn1} are satisfied, then $(M,\pi,N,\ttf)$ is a Poisson quasi-Nijenhuis manifold.
\end{enumerate}
\end{prop}

\begin{proof}
\begin{enumerate}
\item 
Fact \ref{facts} implies  Eq. \eqref{eq:pqn2}. 
By Proposition \ref{prop:qN}, $\big((T^*M)_{\bivector},d_N,\phi\big)$ is
 a quasi-Lie bialgebroid. It is simple to see that its induced bivector field on $M$ as in Proposition 4.8 of \cite{math.DG/0507396} is $\pi_N$. From Proposition 4.8 of \cite{math.DG/0507396}, it follows that $\lie{\pi_N}{\pi_N}=2\pi\diese(\phi)$.
\item 
Since $\lie{\pi}{\pi_N}=0$, Lemma \ref{lem:PN} implies that $\pi\diese\rond d_N=\lie{\pi_N}{\pi\diese(\cdot)}$ and 
Lemma \ref{lem:3} implies that $d_N$ is a derivation of $\lie{\cdot}{\cdot}_{\pi\diese}$. Hence $\pi$ and $N$ are compatible by Lemma \ref{lem:4}.
Since $\pi$ is non-degenerate, we may apply $(\pi\diese)\inv$ to Eq. \eqref{eq:pqn1}. Then, making use of Lemma \ref{lem:PN}, we get back to $d_N^2=\lie{\ttf}{\cdot}_{\pi\diese}$. 
Eq. \eqref{eq:pqn1} and the graded Jacobi identity yield $\lie{\pi_N}{\pi\diese(\phi)}=0$.
Applying $(\pi\diese)\inv$, we get $d_N\ttf=0$.
\end{enumerate}
\end{proof}
 
\begin{cor} \label{prop:5.4}
Let $\omega$ be a symplectic $2$-form  and $\phi$ a closed
$3$-form on $M$. Then $(M,\omega,N,\ttf)$ 
is a symplectic quasi-Nijenhuis manifold if and only if
\begin{gather*}
[\omega_N,\omega_N]=2 \ttf \quad\text{and}\quad d\omega_N=0, \\
\end{gather*}
where $[\cdot,\cdot]$ stands for the Schouten bracket on $\Omega^{\bullet}(M)$ 
induced from the Lie algebroid $(T^*M)_\pi$, and $\omega_N$ is the $2$-form on $M$ 
defined by $$\omega_N(X,Y)=\omega(NX,Y), \quad\forall X,Y\in\XX(M).$$
\end{cor}

\begin{proof}
It is well known that, when $\pi$ is non-degenerate, $\pi\diese$ is an isomorphism 
of differential Gerstenhaber algebras from $(\Omega^{\bullet}(M),d,[\cdot,\cdot])$
to $(\XX^{\bullet}(M),[\pi,\cdot ],[\cdot,\cdot])$ \cite{MR1675117,MR1362125}. 
The conclusion thus follows immediately
from Proposition \ref{prop:3.10} since $\pi\diese\omega_N=\pi_N$.
\end{proof}

\begin{rmk}
Poisson Nijenhuis structures arise naturally in the study of integrable systems.
It would be interesting to find applications of Poisson quasi-Nijenhuis
  structures  in integrable systems as well.
\end{rmk}

\section{Universal lifting theorem}

In this section, we recall the universal lifting theorem and its basic ingredients, as it plays a crucial
role in the following sections. For details, see
\cite{math.DG/0507396}.

Let $\gm\gpd M$ be a Lie groupoid, $A\to M$ its Lie algebroid and 
$\Pi\in \mathfrak X ^k(\gm)$  a $k$-vector field on $\gm$.
Define $F_\Pi \in C^{\infty}(T^\ast \gm\times _\gm \stackrel{(k)}{\ldots }\times _\gm T^\ast \gm)$ by
$$F_\Pi (\mu^1, \ldots , \mu^k)=\Pi (\mu^1, \ldots , \mu^k). $$

\begin{defn}
A $k$-vector field $\Pi \in \mathfrak X ^k (\gm)$ is multiplicative if, and only if,
$F_\Pi$ is a 1-cocycle with respect to
 the groupoid $T^\ast \gm\times _\gm \stackrel{(k)}{\ldots }\times _\gm T^\ast
\gm \gpd A^\ast \times _M \stackrel{(k)}{\ldots }\times _M A^\ast $.
\end{defn}

\begin{rmk}
It is simple to see that a bivector field 
$\Pi$ is multiplicative if, and only if,
 the graph of the multiplication
$\Lambda \subset \gm \times \gm \times\gm$ is coisotropic with respect to
$\Pi\oplus \Pi \oplus \bar{\Pi}$, where
$\bar{\Pi}$ denotes the opposite bivector field to $\Pi$.
\end{rmk}

\begin{ex}
If $P\in \gm (\wedge^k A)$, then $\Vec{P}-\ceV{P}$ is multiplicative,
where $\Vec{P}$ and $\ceV{P}$ denote, respectively, the right and left invariant
$k$-vector fields on $\gm$ corresponding to $P$.
\end{ex}

By ${\mathfrak X}^k_{\mult}(\gm)$ we denote the space of all multiplicative
k-vector fields on $\gm$. And ${\mathfrak X}_{\mult}(\gm)=
\bigoplus_k {\mathfrak X}^k_{\mult}(\gm)$.

\begin{prop}[\cite{math.DG/0507396}]
The vector space ${\mathfrak X}_{\mult}(\gm)$ is closed under the Schouten bracket,
and therefore is a graded Lie algebra.
\end{prop}

It is simple to show that for any given $\Pi\in{\mathfrak X}^k_{\mult}(\gm)$
and any $X\in \gm (\wedge^i A)$,  the $(k+i-1)$-vector field
$[\ceV{X}, \Pi]$ is always left invariant.
Define  $\ceV{\delta_\Pi X} \in  \gm  (\wedge^{(k+i-1)} A)$ by 
$$\ceV{\delta_\Pi X}=[\ceV{X}, \Pi ].$$
Thus one obtains a linear operator $\delta_\Pi : \gm  (\wedge^{i} A)
\to \gm  (\wedge^{(k+i-1)} A)$. Here we use  the following
convention: $\gm  (\wedge^{0} A) \cong C^{\infty}(M)$
and  for any $f\in C^{\infty}(M)$, $\ceV{f}= \beta^* f$.
One easily checks that the following identities are
satisfied
\begin{gather*}
\delta_\Pi (P\wedge Q)=(\delta_\Pi P)\wedge Q+(-1)^{p(k-1)}P\wedge \delta_\Pi Q , \\ 
\delta_\Pi [ P, Q ]=[ \delta_\Pi P, Q ] +(-1)^{(p-1)(k-1)}[P,\delta_\Pi Q] , 
\end{gather*}
for all $P\in \gm (\wedge ^pA)$ and $Q\in\gm (\wedge ^q A)$. This leads
to the following definition of $k$-differentials. 

Recall that for any Lie algebroid $A\to M$, $(\gm (\wedge^{\bullet}A ), \wedge,
[\cdot, \cdot ])$ is a Gerstenhaber algebra \cite{MR1675117}.

\begin{defn}
A k-differential on a Lie algebroid $A$ is a degree $(k-1)$ derivation
of the Gerstenhaber algebra $(\gm(\wedge^{\bullet}A),\wedge,
[\cdot, \cdot ])$. I.e. a linear operator
$$\delta:\sections{\wedge^{\bullet}A}\to\sections{\wedge^{\bullet+(k-1)}A}$$
satisfying
\begin{gather*}
\delta (P\wedge Q)=(\delta P)\wedge Q+(-1)^{p(k-1)}P\wedge \delta Q , \\
\delta [ P, Q ] =[ \delta P, Q  ]+(-1)^{(p-1)(k-1)}[P,\delta Q ] ,
\end{gather*}
for all $P\in\sections{\wedge^p A}$ and $Q\in\sections{\wedge^q A}$.
The set of $k$-differentials on $A$ is denoted by $\cala^k(A)$.
\end{defn}

The space of all multi-differentials $\cala(A)=\bigoplus_k\cala^k(A)$
becomes a graded Lie algebra when endowed with the graded commutator:
$$[\delta_1, \delta_2]=\delta_1\smalcirc \delta_2
-(-1)^{(k-1)(l-1)}\delta_2\smalcirc \delta_1, \qquad\text{where }\delta_1\in \cala^k(A)\text{ and }\delta_2\in \cala^l(A) .$$

Below is a list of basic examples.

\begin{exs}
\begin{enumerate}
\item When $A$ is a Lie algebra $\mathfrak{g}$,
then k-differentials are in one-one correspondence with Lie algebra 1-cocycles $\delta: \mathfrak{g}\to \wedge^k \mathfrak{g}$
with respect to the adjoint action.
\item The $0$-differentials correspond to sections $\phi\in\gm(A^*)$ such that
$d_A\phi=0$, i.e. Lie algebroid 1-cocycles with trivial coefficients.
\item The $1$-differentials correspond to the infinitesimals of Lie algebroid automorphisms.
\item If $P\in\sections{\wedge^k A}$, then $\ad_P=\lie{P}{\cdot}$ is clearly a $k$-differential, which is called the \emph{coboundary} $k$-differential associated to $P$.
\item A Lie bialgebroid can be seen as a Lie algebroid together with a 2-differential of square zero. The converse is also true.
\end{enumerate}
\end{exs}

>From the previous discussion,
we know that there exists a linear map  
$$\mathfrak{X}_{\mult}^{\bullet}(\gm)\to\cala^{\bullet}(A):\Pi\mapsto\delta_{\Pi} ,$$
which is a Lie algebra homomorphism since the graded Jacobi identity satisfied by the Schouten bracket implies that 
\begin{equation} 
[\delta_\Pi, \delta_{\Pi'}]= \delta_{[\Pi, \Pi']} .
\label{eq:lah} \end{equation}
Moreover, one has the following 

\theoremstyle{plain}
\newtheorem*{ULT}{Universal Lifting Theorem}
\begin{ULT}[\cite{math.DG/0507396}]
Assume that $\gm\toto M$ is a target-connected 
and  target-simply connected Lie groupoid with
Lie algebroid $A$. Then $$\mathfrak{X}_{\mult}^{\bullet}(\gm)\to\cala^{\bullet}(A):\Pi\mapsto\delta_{\Pi}$$ 
is an \emph{isomorphism} of graded Lie algebras.
\end{ULT}

\section{Symplectic Nijenhuis groupoids}

\begin{defn}
A symplectic Nijenhuis groupoid is a symplectic groupoid
$(\gm \toto M, \tomega)$ equipped with a  multiplicative
$(1, 1)$-tensor $\tN: T\gm \to T\gm$ such that
$(\gm, \tomega, \tN )$  is a  symplectic Nijenhuis structure.
\end{defn}

The main result of this section is the following

\begin{them}
\label{thm:sym-nij}
\begin{rome}
\item The unit space of a symplectic Nijenhuis groupoid is a Poisson Nijenhuis manifold.
\item Every integrable Poisson Nijenhuis manifold is the unit space of a unique target-connected,  target-simply connected symplectic Nijenhuis groupoid.
\end{rome}
\end{them}

Here, by an integrable Poisson Nijenhuis manifold, we mean the corresponding
 Poisson structure is integrable, i.e. it admits an associated symplectic
groupoid. See \cite{MR1973056, MR2128714} for the solution of the integrability problem
for Poisson manifolds and, more generally, Lie algebroids.

Recall that a Poisson Nijenhuis manifold $(M,\pi,N)$ 
gives rise to a Lie bialgebroid $((T^*M)_\pi,d_N)$
according to Corollary \ref{cor:yvette}.
The following lemma gives a useful characterization of those
Lie bialgebroids arising from Poisson Nijenhuis structures.

\begin{lem}
\label{lem:pnij}
Let $(M,\pi)$ be a Poisson manifold.
A Lie bialgebroid $((T^*M)_\pi,\delta)$
is induced by a Poisson Nijenhuis structure if and only if 
$\lie{\delta}{d}=0$,
where $d$ stands for the de Rham differential.
\end{lem}

\begin{proof}
If $(M,\pi,N)$ is  a Poisson Nijenhuis manifold,
then $d_N =i_N\smalcirc d-d\smalcirc i_N$. Thus 
$$[d_N,d]=d_N\smalcirc d+d\smalcirc d_N
=(i_N\smalcirc d-d\smalcirc i_N)\smalcirc d+d\smalcirc(i_N\smalcirc d-d
\smalcirc i_N)=0.$$

Conversely, given a Lie bialgebroid $((T^*M)_\pi,\delta)$
such that $\lie{\delta}{d}=0$, one obtains
a Lie algebroid structure on $TM$.
Let $N:TM\to TM$ be its anchor map.
Thus $\delta=d_N:C^\infty(M)\to\Omega^1(M)$.
Since $\lie{\delta}{d}=0$, we have $\forall f\in C^\infty(M)$,
$\delta(df)=-d\delta f=-d d_Nf=d_N(df)$. It thus follows
that $\delta=d_N$ on any differential forms
since both $\delta$ and $d_N$ are derivations and they agree 
on 0- and exact 1-forms.
According to Corollary \ref{cor:yvette}, it follows that
$(M,\pi,N)$ is a Poisson Nijenhuis manifold.
\end{proof}

\begin{proof}[Proof of Theorem \ref{thm:sym-nij}]
\begin{rome}
\item \textsc{From symplectic Nijenhuis groupoids to Poisson Nijenhuis manifolds.}
Assume that $(\gm,\tomega,\tN)$ is a symplectic Nijenhuis
groupoid. Let $\tpi$ be the bivector field on $\gm$ which
is the inverse of $\tomega$ and $\tpi_\tN\in\XX^2(\gm)$ be the bivector
field defined by $\tpi_\tN\diese=\tN  \smalcirc\tpi\diese$.
\begin{itemize}
\item Since $\lie{\tpi}{\tpi}=0$, the induced bivector field $\pi=t_*\tpi$ on the base manifold 
of the symplectic groupoid $\gm\toto M$ is Poisson \cite{MR866024}. The 
Lie algebroid of $\gm\to M$ is isomorphic to 
$(T^*M)_\pi$ \cite{MR996653}. And the multiplicative bivector field
$\tpi$ corresponds to a $2$-differential on $(T^*M)_\pi$, which is
the de Rham differential $d$. 
That is, $((T^*M)_\pi,d)$ is the Lie bialgebroid corresponding to
the symplectic groupoid $(\gm,\tomega)$.
\item As pointed out in Fact \ref{facts}, $\tpi_\tN$ 
is a Poisson tensor on $\gm$ \cite{MR773513,MR1077465,MR1390832}.  
Moreover, $\tpi_\tN$ is a multiplicative bivector field since 
$\tN$ is a multiplicative $(1,1)$-tensor and
$\tpi$ is a multiplicative bivector field.
In other words, $(\gm,\tilde{\pi}_\tN)$ is a Poisson groupoid \cite{MR1262213}.
Let $\delta_{\tpi_\tN}:\Omega^{\bullet}(M)\to\Omega^{\bullet+1}(M)$ be the
$2$-differential on $(T^*M)_\pi$ induced by the 
multiplicative Poisson bivector field $\tpi_\tN$ on $\gm$.
Since $\lie{\tpi_\tN}{\tpi_\tN}=0$, the universal lifting theorem implies that 
$$0=\delta_{\lie{\tpi_\tN}{\tpi_\tN}}=\lie{\delta_{\tpi_\tN}}{\delta_{\tpi_\tN}}
=\delta_{\tpi_\tN}\smalcirc \delta_{\tpi_\tN}+\delta_{\tpi_\tN}\smalcirc
\delta_{\tpi_\tN}=2\delta_{\tpi_\tN}^2 .$$
Thus, $((T^*M)_\pi,\delta_{\tpi_\tN})$ is a Lie bialgebroid. 
\item Likewise, it is standard that $\lie{\tpi_\tN}{\tpi}=0$. Thus the universal 
lifting theorem implies that $\lie{\delta_{\tpi_\tN}}{d}=0$. 
According to Lemma \ref{lem:pnij}, $\delta_{\tpi_\tN}=d_N$ for some 
Nijenhuis tensor $N$ on $M$ and $(M,\pi,N)$ is a Poisson Nijenhuis manifold.
\end{itemize}
\item \textsc{From Poisson Nijenhuis manifolds to symplectic Nijenhuis groupoids.}
Given a Poisson Nijenhuis manifold
$(M,\pi,N)$, then $((T^*M)_\pi,d_N)$ is a Lie
bialgebroid by Corollary \ref{cor:yvette}. Assume that $(T^*M)_{\pi}$ is integrable  (see \cite{MR1973056, MR2128714} for the integrability condition)
and $(\gm\toto M,\tomega)$ is a target-connected
and target simply-connected symplectic groupoid of $M$.
Since $d_N^2=0$ and $\lie{d_N}{d}=0$, the universal lifting theorem implies that 
$d_N$ corresponds to a multiplicative Poisson bivector field $\tpi_\tN$ on $\gm$ such that 
$\lie{\tpi_\tN}{\tpi}=0$,
where $\tpi$ is the Poisson tensor on $\gm$ inverse to $\tomega$. 
Let $\tN=\tpi_\tN\diese\smalcirc\tomega\bemol:T\gm\to T\gm$.
Then it is clear that $\tN$ is a multiplicative
$(1,1)$-tensor, and $(\gm,\tomega,\tN)$ is a symplectic Nijenhuis
groupoid.
\end{rome}
Since these two constructions are inverse to each other, the theorem is proved.
\end{proof}

\section{Symplectic quasi-Nijenhuis groupoids}

The goal of this section is to generalize
Theorem \ref{thm:sym-nij} to the quasi-setting. More precisely,
we will give an integration theorem for
Poisson quasi-Nijenhuis manifolds.

\begin{defn}
A symplectic quasi-Nijenhuis groupoid is a symplectic groupoid
$(\gm\toto M,\tomega)$ equipped with a multiplicative $(1,1)$-tensor 
$\tN:T\gm\to T\gm$ and a closed $3$-form $\ttf\in\Omega^3(M)$ 
such that $\big(\gm,\tomega,\tN,t^*\ttf-s^*\ttf\big)$ is a symplectic 
quasi-Nijenhuis structure.
\end{defn}

The following result is a generalization of Theorem \ref{thm:sym-nij}.

\begin{them} \label{thm:sym-nij1}
\begin{rome}
\item The unit space of a symplectic quasi-Nijenhuis groupoid
is a Poisson quasi-Nijenhuis manifold. 
\item Every integrable Poisson quasi-Nijenhuis manifold $(M,\pi,N,\ttf)$ is the unit space of a 
unique target-connected and target-simply connected symplectic quasi-Nijenhuis groupoid
$\big(\gm\toto M,\tomega,\tN,t^*\ttf-s^*\ttf\big)$.
\end{rome}
\end{them}
\begin{proof}
The proof is similar to that of Theorem \ref{thm:sym-nij},
so we will merely sketch it.

Assume that  
$(M,\pi,N,\ttf)$ is an integrable Poisson quasi-Nijenhuis manifold.
Let $\gm\toto M$ be a target-connected and  target-simply
connected groupoid integrating the Lie algebroid
$(T^*M)_\pi$. By Proposition \ref{prop:qN}, 
$\big((T^*M)_\pi,d_N,\ttf\big)$ is a quasi-Lie bialgebroid, 
which integrates to a quasi-Poisson groupoid by the universal lifting theorem. 
Let $\tpi_\tN\in\XX(\gm)$ be the bivector 
field on $\gm$ corresponding to $d_N$. Then we have
\[ \thalf \lie{\tpi_\tN}{\tpi_\tN}=\Vec{\ttf}-\ceV{\ttf} .\]
On the other hand, we know that $\gm\toto M$ 
is a symplectic groupoid, whose corresponding
Lie bialgebroid is $\big((T^*M)_\pi,d\big)$. 
The symplectic form on $\gm$ is denoted by $\tomega$. 
Let $\tpi\in\XX^2(\gm)$ be its corresponding Poisson tensor. 
Since $\lie{d_N}{d}=0$, we have $\lie{\tpi_\tN}{\tpi}=0$ according 
to the universal lifting theorem. Let 
$\tN=\tpi_\tN\diese\smalcirc\tomega\bemol:T\gm\to T\gm$.
Then it is clear that $\tN$ is a multiplicative
$(1,1)$-tensor. Since $\Vec{\ttf}-\ceV{\ttf}=\tpi\diese
(t^*\ttf-s^*\ttf)$, from Proposition \ref{prop:3.10},
it follows that $\big(\gm,\tomega,\tN,t^*\ttf-s^*\ttf\big)$ 
is a symplectic quasi-Nijenhuis groupoid.

The other direction can be proved by going backwards. 
\end{proof}

\begin{rmk}
Note that  $\tomega\bemol (\tpi_\tN )$ is a multiplicative $2$-form
on $\gm\toto M$. It would be interesting to see what is the corresponding
Dirac structure on $M$ and how the integration result in \cite{MR2068969}
can be applied to this situation.
\end{rmk}

\section{Generalized complex structures} \label{sec:GCS}

This section is devoted to the investigation of the
relationship between generalized complex structures
and Poisson quasi-Nijenhuis structures.
Let us first recall the definition of generalized complex structures \cite{MR2013140,math.DG/0401221}. 

\begin{defn}
A generalized complex structure on a manifold $M$ is a bundle map 
$$J:TM\oplus T^*M\to TM\oplus T^*M$$ satisfying the algebraic properties 
\begin{equation} J^2=-I \qquad\text{and}\qquad \ip{Jv}{Jw}=\ip{v}{w} \label{eq:J1} \end{equation}
and the integrability condition 
\[ \std{Jv}{Jw}-\std{v}{w}-J(\std{Jv}{w}+\std{v}{Jw})=0 \]
$\forall v,w\in\Gamma(TM\oplus T^*M)$. Here $\ip{\cdot}{\cdot}$
and $\std{\cdot}{\cdot}$ are the pairing and bracket
on the standard Courant algebroid $TM\oplus T^*M$ as
in Example \ref{ex:1.2}.
\end{defn}

The first two algebraic  conditions \eqref{eq:J1} imply that $J$
 must be of the form 
\begin{equation} \label{eq:J}
J=\begin{pmatrix}N & \pi^{\sharp} \\ \sigma_{\flat} & -N^*\end{pmatrix}
\end{equation} 
where $\pi\in\mathfrak{X}^2(M)$ is a bivector field, $\sigma\in\Omega^2(M)$
is a 2-form and $N:TM\to TM$ is a $(1,1)$-tensor.
Here $\sigma_{\flat}:TM\to T^*M$ is the map given by
$(\sigma_{\flat}X)(Y)=\sigma(X,Y)$, $\forall X,Y\in\XX(M)$.

On the other hand, a Courant algebroid can be deformed using
a bundle map $J$. More precisely, 
let $(E,\ip{\cdot}{\cdot},\courant{\cdot}{\cdot},\anchor)$ be a Courant algebroid over $M$ and let 
$$\xymatrix{ E \ar[r]^J \ar[d] & E \ar[d] \\ M \ar[r]_{\idn} & M }$$ be a vector bundle automorphism of $E\to M$.
Consider
\begin{itemize}
\item the inner product 
\begin{equation*}
\label{eq:N1}
\ipJ{A}{B}=\ip{JA}{JB},
\end{equation*}
\item the bracket
 \begin{equation}
\label{eq:N2}
\courant{A}{B}_J=\courant{JA}{B}+\courant{A}{JB}-J\courant{A}{B}
\end{equation}
\item and the bundle map 
\begin{equation*}
\label{eq:N3}
\anchor_J=\anchor\rond J
\end{equation*}
\end{itemize}
induced by $J$.

A natural question is

\begin{question} \label{quest}
When is the quadruple $(E,\ipJ{\cdot}{\cdot},\courant{\cdot}{\cdot}_J,\anchorJ)$ 
still a Courant algebroid?
\end{question}

The next proposition gives a trivial sufficient condition.

\begin{prop} \label{prop:2.8}
The quadruple $(E,\ip{\cdot}{\cdot}_J,\courant{\cdot}{\cdot}_J,\anchor_J)$ is a Courant algebroid if
$$\courant{JA}{JB}+J^2\courant{A}{B}-J\big(\courant{JA}{B}+\courant{A}{JB}\big)=0, \quad \forall A,B\in\sections{E}.$$
Moreover, in this case, $J$ is a Courant algebroid isomorphism from 
$(E,\ip{\cdot}{\cdot}_J,\courant{\cdot}{\cdot}_J,\anchor_J)$ 
to $(E,\ip{\cdot}{\cdot},\courant{\cdot}{\cdot},\anchor)$.
\end{prop}

We now give an answer to Question~\ref{quest} in the special case
 of the standard Courant algebroid $TM\oplus T^*M$, 
where  $J$ satisfies Eqs. \eqref{eq:J1}, and is given by Eq. \eqref{eq:J} .

\begin{lem}
\label{lem:3.3}
Assume that $J:TM\oplus T^*M\to TM\oplus T^*M$ is given by Eq. \eqref{eq:J}.
Let $\std{\cdot}{\cdot}_J$ be the deformed bracket on $\XX(M)\oplus \Omega^1 (M)$
as in Eq. \eqref{eq:N2}.
Then, for all $X,Y\in\XX(M)$ and $\xi,\eta\in\Omega^1(M)$, we have
\begin{gather}
\std{\xi}{\eta}_J=[\xi,\eta]_{\pi\diese} \label{eq:xieta}\\
\std{X}{Y}_J=[X,Y]_N +(d\sigma)(X,Y,\cdot) \label{eq:xy}\\
\std{X}{\xi}_J=\big([X,\pi\diese\xi]-\pi\diese(\derlie{X}\xi-\thalf d(\xi X))\big)
+\big(\derlie{NX}\xi-\derlie{X}(N^T\xi)+N^T(\derlie{X}\xi-\thalf d(\xi X))\big) \label{eq:Xxi}
\end{gather}
\end{lem}

\begin{proof}
This follows from a straightforward computation using
Eqs. \eqref{eq:bracket} and \eqref{eq:N2}, and is left for the reader.
\end{proof}

\begin{prop} \label{prop:3.4}
Let $J:TM\oplus T^*M\to TM\oplus T^*M$
be a bundle map which satisfies Eqs. \eqref{eq:J1}, 
and is given by Eq. \eqref{eq:J}.
Then $\find$ is a Courant algebroid if, and only if,
$(M,\pi,N,d\sigma)$ is a Poisson quasi-Nijenhuis manifold.
And in this case, $\find$ is  naturally
identified with the double of the quasi-Lie bialgebroid 
$\big((T^*M)_{\pi},d_N,d\sigma\big)$.
\end{prop}

\begin{proof}
Assume that $(TM\oplus T^*M,\ip{\cdot}{\cdot}_J,\std{\cdot}{\cdot}_J,\anchor_J)$ is a Courant
 algebroid. It is clear that $A:=T^*M$ and $B: =TM$ are transversal, maximal 
isotropic subbundles.  By Eq. \eqref{eq:xieta}, $A=T^*M$ is a Dirac structure with the induced bracket $\lie{\cdot}{\cdot}_{\bivector\diese}$. Thus, according to Theorem \ref{thm:1.6}, we obtain a quasi-Lie bialgebroid. The construction of the corresponding derivation $\delta$ of $(\Omega^{\bullet}(M),\wedge,\lie{\cdot}{\cdot}_{\bivector\diese})$ and twisting 3-form $\phi$ was outlined in the proof of Theorem \ref{thm:1.6}. In the present situation, we have $$\rho_B(X)=\rho_J(X)=\rho(JX)=\rho(NX+\sigma\bemol X)=NX, \quad \forall X\in TM$$ and, combining Eqs. \eqref{eq:xy} and \eqref{eq:5bis}, $$\lie{X}{Y}_B=\lie{X}{Y}_N, \quad \forall X,Y\in\vf(M) .$$
Therefore, comparing Eqs. \eqref{eq:deltaA} and \eqref{eq:8bis}, we conclude that $\delta=d_N$.
And, combining Eqs. \eqref{eq:5ter} and \eqref{eq:xy}, we get 
\begin{multline*} 
\phi(X,Y,Z)=2\ip{\std{X}{Y}_J}{Z}_J=2\ip{J\std{X}{Y}_J}{JZ}=2\ip{\std{X}{Y}_J}{Z} \\ 
=2\ip{\lie{X}{Y}_N+d\sigma(X,Y,\cdot)}{Z}=d\sigma(X,Y,Z), \quad\forall X,Y,Z\in\vf(M).
\end{multline*}
Hence $\big((T^*M)_{\bivector}, d_N, d\sigma\big)$ is a quasi-Lie bialgebroid or, equivalently according to Proposition \ref{prop:qN}, $(M,\bivector,N,d\sigma)$ is a Poisson quasi-Nijenhuis manifold.

Conversely, assume that $(M,\pi,N,d\sigma)$ is a Poisson quasi-Nijenhuis manifold.
By Proposition \ref{prop:qN}, $((T^*M)_\pi,d_N,d\sigma)$ is a quasi-Lie bialgebroid. 
Its double $E$ is a Courant algebroid. We will show that $E$ 
is indeed isomorphic to $\find$. First, it is simple to check that 
their anchors and non-degenerate symmetric pairings coincide. 
It remains to check that their brackets coincide. According to Eq.
\eqref{eq:5quater}, the bracket $\courant{\cdot}{\cdot}$ on $\sections{E}$ is given by
\begin{gather}
\courant{\xi}{\eta}=[\xi,\eta]_\pi \label{eq:xieta1} \\
\courant{X}{Y}=\lie{X}{Y}_N+(d\sigma)(X,Y,\cdot) \label{eq:xy1} \\
\courant{X}{\xi}=\big(i_X\delta_{TM}\xi+\thalf\delta_{TM}(\xi X)\big)
-\big(i_{\xi}\delta_{T^*M} X+\thalf\delta_{T^*M}(\xi X)\big) \label{eq:Xxi1}
\end{gather}
for all $X,Y\in\XX(M)$ and $\xi,\eta\in\Omega^1(M)$.
In our case, we have
$$\delta_{T^*M}=\lie{\pi}{\cdot} \qquad\text{and}\qquad \delta_{TM}=d_N.$$
It follows from a straightforward verification that the right hand sides 
of Eqs. \eqref{eq:xieta}-\eqref{eq:Xxi} and  
\eqref{eq:xieta1}-\eqref{eq:Xxi1} coincide.
Therefore, $\find$ is indeed a Courant algebroid.
\end{proof}

We are now ready to state the main result of this section.

\begin{them}
Assume that $J:TM\oplus T^*M\to TM\oplus T^*M$ as given by Eq. \eqref{eq:J}
satisfies Eqs \eqref{eq:J1}. Then the following are equivalent
\begin{itemize}
\item $J$ is a generalized complex structure;
\item $(M,\pi,N,d\sigma)$ is a Poisson quasi-Nijenhuis manifold such that
$$(TM)_N\oplus (T^*M)_{\pi}\xrightarrow{J}TM\oplus T^*M$$
is a Courant algebroid isomorphism. 
\end{itemize}
Here $(TM)_N\oplus (T^*M)_{\pi}$ denotes 
the Courant algebroid corresponding to the
quasi-Lie bialgebroid $((T^*M)_\pi,d_N,d\sigma)$.
\end{them}

\begin{proof}
By Proposition \ref{prop:2.8}, $J$ is a generalized complex structure
 if, and only if, $\find$ is a Courant algebroid and
 $(TM\oplus T^*M,\ip{\cdot}{\cdot}_J,\std{}{}_J,\rho_J)\xrightarrow{J}(TM\oplus T^*M,\ip{\cdot}{\cdot},\std{\cdot}{\cdot},\anchor)$ 
is a Courant algebroid isomorphism.
The result follows  immediately from Proposition \ref{prop:3.4}.
\end{proof}

Since any generalized complex structure naturally
gives rise to a Poisson quasi-Nijenhuis
manifold, as an immediate consequence of Theorem~\ref{thm:sym-nij1}, we have
the following

\begin{them}
\label{thm:J}
Let $J$ be a generalized complex structure as given by
Eq. \eqref{eq:J}, and $(\gm\toto M,\tilde{\omega})$ a target-connected and
 target-simply connected symplectic groupoid integrating
$(T^*M)_\pi$. Then there is a multiplicative
$(1,1)$-tensor $\tN$ on $\gm$ such that 
$(\gm\toto M,\tilde{\omega}, \tN,t^*d\sigma-s^*d\sigma)$ is a 
symplectic quasi-Nijenhuis groupoid.
\end{them}

\begin{rmk}
Note that Theorem 3.3-3.4 in \cite{math.DG/0412097} essentially
imply our Theorem \ref{thm:J}.
Our proof is conceptual, while Crainic used a direct argument. 
It would be interesting to see how Theorem 3.4 (ii) in \cite{math.DG/0412097} 
can be proved conceptually.
\end{rmk}

\bibliographystyle{hamsplain}



\def\cprime{$'$}
\providecommand{\bysame}{\leavevmode\hbox to3em{\hrulefill}\thinspace}
\providecommand{\href}[2]{#2}

\end{document}